\documentclass{amsart}

\usepackage{xypic}
\input xy
\xyoption{all}
\usepackage{epsfig}
\usepackage{amsthm}
\usepackage{amssymb}
\usepackage{amsmath}
\usepackage{amscd}
\usepackage{color}

\newcommand{\bg}{\begin{equation}}
\newcommand{\ed}{\end{equation}}
\newcommand{\bga}{\begin{eqnarray}}
\newcommand{\eda}{\end{eqnarray}}
\newcommand{\pf}{\textbf{Proof:\ }}

\def\cbdu{\par{\raggedleft$\Box$\par}}

\newtheorem {Theorem}  {Theorem}

\numberwithin{Theorem}{section}

\newtheorem {Lemma}[Theorem]  {Lemma}

\theoremstyle{definition}

\theoremstyle{remark}
\newtheorem{Remark}[Theorem]{\bf Remark}

\expandafter\chardef\csname pre amssym.def
at\endcsname=\the\catcode`\@ \catcode`\@=11
\def\undefine#1{\let#1\undefined}
\def\newsymbol#1#2#3#4#5{\let\next@\relax
 \ifnum#2=\@ne\let\next@\msafam@\else
 \ifnum#2=\tw@\let\next@\msbfam@\fi\fi
 \mathchardef#1="#3\next@#4#5}
\def\mathhexbox@#1#2#3{\relax
 \ifmmode\mathpalette{}{\m@th\mathchar"#1#2#3}%
 \else\leavevmode\hbox{$\m@th\mathchar"#1#2#3$}\fi}
\def\hexnumber@#1{\ifcase#1 0\or 1\or 2\or 3\or 4\or 5\or 6\or 7\or 8\or
 9\or A\or B\or C\or D\or E\or F\fi}

\font\teneufm=eufm10 \font\seveneufm=eufm7 \font\fiveeufm=eufm5
\newfam\eufmfam
\textfont\eufmfam=\teneufm \scriptfont\eufmfam=\seveneufm
\scriptscriptfont\eufmfam=\fiveeufm

\catcode`\@=\csname pre amssym.def at\endcsname

\newcounter{remark}
\newenvironment{remark}
{\medskip \stepcounter{remark} \noindent \textit{Remark
\arabic{section}.\arabic{remark}.}}{\rm \cbdu}

\newcommand{\supp}{{\mathit supp}\,}

\renewcommand{\d}{\delta}

\renewcommand{\l}{\lambda}

\newcommand{\R}{\mathbf{R}}

\def  \R   {{\mathbb R}}

\def  \12  {{\frac{1}{2}}}

\def  \l   {\langle}

\def\build#1_#2^#3{\mathrel{\mathop{\kern 0pt#1}\limits_{#2}^{#3}}}

\begin{document}

\title[Ill-posedness for NSE and MHD ]{Ill-posedness of  the Navier-Stokes and magneto-hydrodynamics systems}


\author [Alexey Cheskidov]{Alexey Cheskidov}
\address{Department of Mathematics, University of Illinois, Chicago, IL 60607,USA}
\email{acheskid@uic.edu}

\author [Mimi Dai]{Mimi Dai}
\address{Department of Mathematics, University of Illinois, Chicago, IL 60607,USA}
\email{mdai@uic.edu}

\thanks{The work of Alexey Cheskidov was partially supported by NSF Grant DMS-1108864}





\begin{abstract}
We demonstrate that the three
dimensional incompressible magneto-hydrodynamics (MHD) system is ill-posed due to the discontinuity of weak solutions in a wide range of spaces.
Specifically, we construct initial data which has finite energy and is small in certain spaces, such that any Leray-Hopf type of weak solution to the MHD system starting from this initial data is discontinuous at time $t=0$ in such spaces. An analogous result is also obtained for the Navier-Stokes equation which extends the previous result of ill-posedness in $\dot B^{-1}_{\infty,\infty}$ by Cheskidov and Shvydkoy to spaces that are not necessarily critical. The region of the spaces where the norm inflation occurs almost touches $L^2$.

\bigskip


\hspace{0.02cm}CLASSIFICATION CODE: 76D03, 35Q35.
\end{abstract}

\maketitle

\section{Introduction}

The three dimensional incompressible
magneto-hydrodynamics (MHD) system is given by:
\begin{equation}\label{MHD}
\begin{split}
u_t- \mu\triangle u +u\cdot\nabla u-b\cdot\nabla b+\nabla p=0,\\
b_t- \nu\triangle b+u\cdot\nabla b-b\cdot\nabla u=0,\\
\nabla \cdot u=0, \ \ \ \nabla \cdot b=0, 
\end{split}
\end{equation}
with the initial conditions
\begin{equation}\label{init}
\begin{split}
u(x,0)=u_0(x),\;\;\; b(x,0)=b_0(x),\\
\nabla\cdot u_0=0, \;\;\; \nabla\cdot b_0=0
\end{split}
\end{equation}
where $x\in \Omega = \mathbb{T}^3$, $t\geq 0$, $u$ is the fluid velocity, $p$ is the pressure of the fluid, and $b$
is the magnetic field. The parameter $\mu$ denotes the kinematic viscosity coefficient of the fluid and $\nu$ denotes the reciprocal of the magnetic Reynolds number. 
When the magnetic field $b(x,t)$ vanishes,
the incompressible MHD system reduces to the incompressible Navier-Stokes
equation (NSE). In the case where the domain $\Omega$ is the whole space, the solutions to the MHD system share the same scaling
property of the solutions to the NSE, that is,
\begin{equation}
u_\lambda(x,t) =\lambda u(\lambda x, \lambda^2t), \ b_\lambda(x,t)
=\lambda b(\lambda x, \lambda^2t), \ p_\lambda(x,t) =\lambda^2
p(\lambda x, \lambda^2t) \notag
\end{equation}
solve (\ref{MHD}) with the initial data
\begin{equation}
u_{0\lambda} =\lambda u_0(\lambda x), \ b_{0\lambda} =\lambda
b_0(\lambda x), \notag
\end{equation}
if $(u(x, t), b(x, t))$ solves (\ref{MHD})
with the initial data $(u_0(x), b_0(x))$. A space that is
invariant under the above scaling is called critical space.
Examples of critical spaces 
associated with the above scaling in three dimension are
$$
\dot{H}^\frac 12 \hookrightarrow L^3\hookrightarrow\dot B^{-1+\frac{3}{p}}_{p|p<\infty,\infty}\hookrightarrow BMO^{-1}
\hookrightarrow \dot{B}^{-1}_{\infty,\infty}.
$$
Notice that $\dot{B}^{-1}_{\infty,\infty}$ is the largest critical space for both the NSE and the MHD system.
In the periodic case there is no distinction between homogeneous and non-homogeneous spaces, so $ B^{-1+\frac{3}{p}}_{p,\infty}$ are also called critical.

The study of the Navier-Stokes equations in
critical spaces has been a focus of the research activity
since the initial work of Kato \cite{Kato}.  In 2001, Koch and Tataru \cite{KT}
established the global well-posedness of the classical Navier-Stokes equations with
small initial data in the space $BMO^{-1}$. Then the question whether this result can be
extended to the largest critical space $\dot{B}^{-1}_{\infty,\infty}$ had become of great interest.
The first indication that such an extension might not be possible came in the work by
Bourgain and Pavlovi\'{c} \cite{BP} who showed the
norm inflation for the classical Navier-Stokes equations in $\dot{B}^{-1}_{\infty,\infty}$.
More precisely, they constructed arbitrarily small initial
data in $\dot{B}^{-1}_{\infty,\infty}$, such that mild solutions with this data become arbitrarily large in $\dot{B}^{-1}_{\infty,\infty}$ after an arbitrarily short time. This result was later extended to generalized Besov spaces smaller than $B^{-1}_{\infty,p}$, $p>2$ by Yoneda \cite{Yo}.
Moreover, in \cite{CS} Cheskidov and Shvydkoy proved the
existence of discontinuous Leray-Hopf solutions of the Navier-Stokes equations in  $\dot{B}^{-1}_{\infty,\infty}$ with arbitrarily small initial data.
Contrary to the Bourgain-Pavlovi\'{c} construction where the energy transfers from high to low modes to produce the norm inflation, the norm discontinuity in \cite{CS}
is due to the forward energy cascade generated by local interactions.
 In \cite{CD} Cheskidov and Dai considered fractional Navier-Stokes equations and showed that the natural space for norm inflation is critical only when the power of the Laplacian is one.
When the power od the laplacian is larger than one, the norm inflation occurs not only in critical spaces, but also in subcritical and supercritical.

For the MHD system, Miao, Yuan and Zhang \cite{MYZ} proved the
existence of a global mild solution in $(BMO^{-1})^2$ for small initial
data and uniqueness of such solution in $\left(C([0,\infty); BMO^{-1})\right)^2$. 
Later, Dai, Qing and Schonbek \cite{DQSM} established several different types of ``norm inflation"  phenomena for the three dimensional MHD system in the largest critical space $(\dot{B}^{-1}_{\infty,\infty})^2$, by adopting the idea of \cite{BP}. Since the MHD system describes the coupling of velocity field and magnetic field, the authors were able to construct different initial data to produce different types of ``norm inflation". In particular, the magnetic field can develop norm inflation in short time even when the velocity remains small and vice versa. In \cite{CDmhd} Cheskidov and Dai used their approach  from
\cite{CD} to extend the norm inflation results to even wider range of spaces that included critical, subcritical, and supercritical.

In this paper we further investigate the method of \cite{CS} to study the ill-posedness problem of the NSE and the MHD system in a large class of spaces which may contain critical, supercritical and subcritical spaces.  First, modifying the initial data construction, we are able to obtain discontinuous weak solutions to the NSE in certain Besov spaces. Namely, we prove that

\begin{Theorem}\label{thm-nse}
Let $1<r\leq\infty$ and $\frac 32<\theta\leq 2$ which satisfy either 
\[\left\{ 2\leq\theta+\frac3r<3\right\}; \quad \mbox { or }\]
\[\left\{ r\geq \frac32, \ \qquad  3\leq\theta+\frac3r<4,  \ \qquad  2\theta+\frac3r\leq \frac{11}2\right\}.\]
There exists an initial data $u_0\in B_{r,\infty}^{\frac{3}{r}+\theta-3}$ (that depends only on $\theta$),  such that every Leray-Hopf weak solution $u\in C_w([0,T); L^2) \cap L^2([0,T); H)$  to the NSE satisfies
\begin{equation}
\limsup_{t \to 0+} \|u(t) - u_0\|_{B^{\frac3r+\theta-3}_{r,\infty}} \geq \delta
\end{equation}
for an absolute constant $\delta$. 
\end{Theorem}

\begin{figure}
\label{fig}
\centering
\includegraphics[width=3.5in, trim =0in 6in 0.5in 1.5in, clip =true]{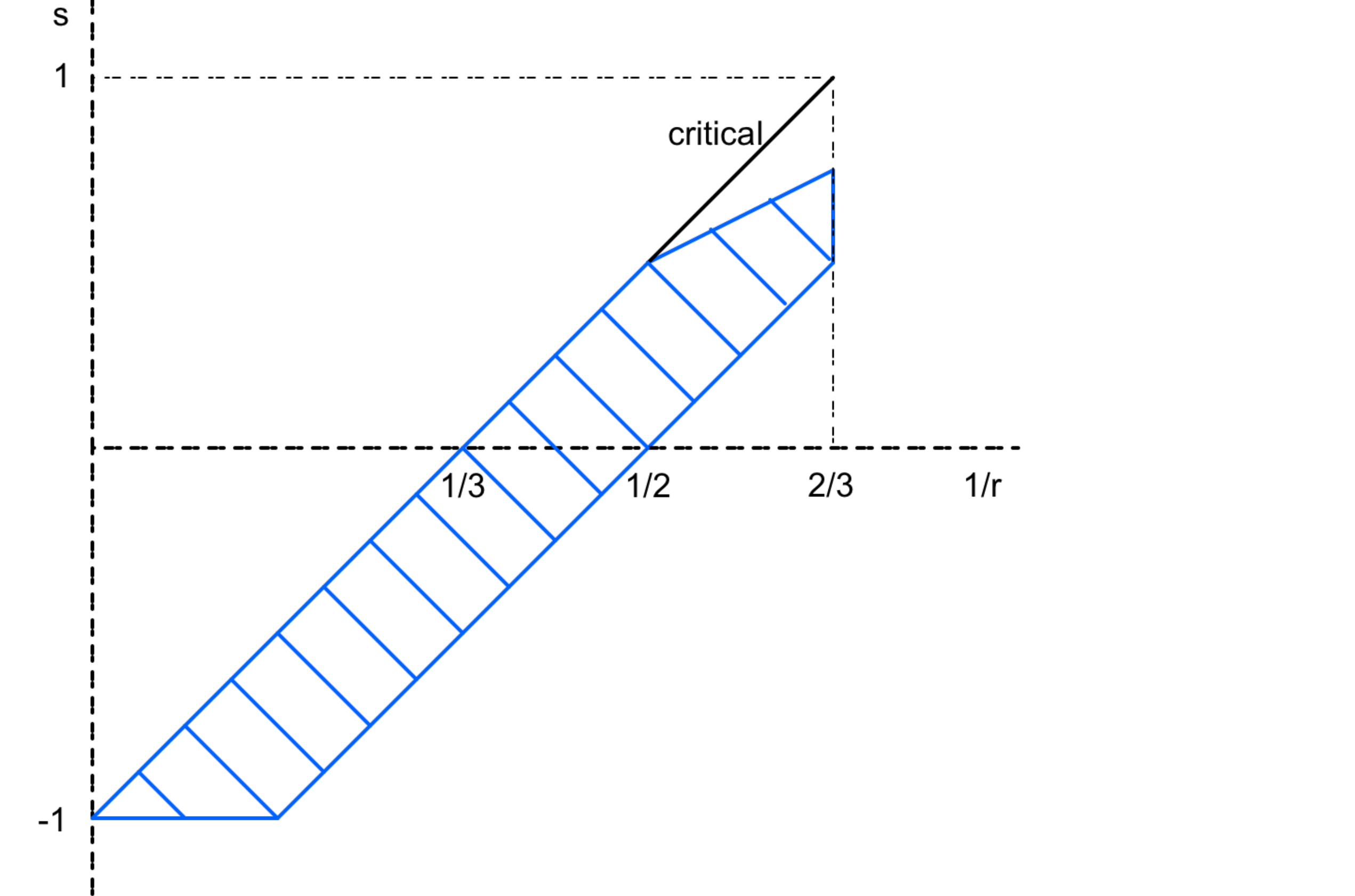}
\caption{The region where the discontinuity for the NSE occurs on the plane of the smoothness index $s$ vs $1/r$.}
\label{fig3}
\end{figure}

The region of the spaces where the discontinuity of the NSE occurs is diagrammed as in Figure \ref{fig}. One can see that the NSE develops discontinuous weak solutions in both critical and supercritical spaces. The region contains the largest critical space $B^{-1}_{\infty,\infty}$, but also $B^0_{2,\infty}$ that has the same scaling as $L^2$.

In the mean time, we obtain  the same type of ill-posedness for the MHD system as follows:

\begin{Theorem}\label{thm} 
Let $1<r\leq\infty$, $\frac 32<\theta\leq 2$ and $\gamma>\frac 32$ with $\theta+\gamma\leq 4$. In addition, the triplet $(r,\theta,\gamma)$ satisfies either
\begin{equation}\notag
\left\{\gamma\leq \frac52,  \quad 2\leq\theta+\frac3r<3,  \quad \gamma+\frac3r<4\right\}; \quad \mbox {or}
\end{equation}
\begin{equation}\notag
\left\{r\geq\frac32, \, \theta+\frac3r\geq3, \, \gamma+\frac3r<4, \, \theta+\gamma+\frac3r\leq\frac{11}2\right\}.
\end{equation}
There exists an initial data $(u_0,b_0)\in B_{r,\infty}^{\frac{3}{r}+\theta-3}\times B_{r,\infty}^{\frac{3}{r}+\gamma-3}$, such that every Leray-Hopf weak solution $(u(t),b(t))$ to (\ref{MHD})-(\ref{init}) satisfies
\begin{equation}\label{ineq:main}
\limsup_{t\to 0^+}\left(\|u(t)-u_0\|_{B^{\frac3r+\theta-3}_{r,\infty}}+\|b(t)-b_0\|_{B^{\frac3r+\gamma-3}_{r,\infty}}\right)\geq\delta
\end{equation}
for an absolute constant $\delta$.
\end{Theorem}


\begin{remark} We point out that the assumption on the parameter triplet $(r,\theta,\gamma)$ may not be optimal. In fact, in the proof of the theorem in Section \ref{sec:MHD}, the assumption guarantees that a ``jump" of $\|b(t)\|_{L^2}$ occurs in the contradiction argument. One can verify that different (complimentary) assumption on $(r,\theta,\gamma)$ may yield a ``jump" of  $\|u(t)\|_{L^2}$. We do not include the alternate assumption in the statement of the theorem due to the complication. However, under the current assumption, one can already see the discontinuity occurs in a wide range of spaces including critical, supercritical and subcritical ones. 
\end{remark}

\medskip

The rest of the paper is organized as: in Section \ref{sec:pre} we introduce some notations that shall be used throughout the paper and some auxiliary results; in Section \ref{sec:data} we present the initial data construction for both the NSE and the MHD system; Section \ref{sec:NSE} provides a brief proof of Theorem \ref{thm-nse}; Section \ref{sec:MHD} is devoted to proving Theorem \ref{thm}.

\bigskip

\section{Preliminaries and auxiliary results}
\label{sec:pre}

\subsection{Notation}
\label{sec:notation}
We denote by $A\lesssim B$ an estimate of the form $A\leq C B$ with
some constant $C$, and by $A\sim B$ an estimate of the form $C_1
B\leq A\leq C_2 B$ with some constants $C_1$, $C_2$. 
We denote $\|\cdot\|_p=\|\cdot\|_{L^p(\mathbb T^n)}$ and the trilinear term
\begin{equation}\label{tril}
\mathcal B(u,v,w)=\int_{\mathbb{T}^3} u\otimes v:\nabla w\, dx=\int_{\mathbb{T}^3}v_i\partial_iw_ju_j\, dx.
\end{equation}

\subsection{Littlewood-Paley decomposition}
\label{sec:LPD}
The techniques presented in this paper rely strongly on the Littlewood-Paley decomposition. We recall the Littlewood-Paley decomposition theory briefly. For a more detailed description on this theory we refer the readers to the books by Bahouri, Chemin and Danchin \cite{BCD} and Grafakos \cite{Gr}. 

We denote $\lambda_q=2^q$ for integers $q$. A nonnegative radial function $\chi\in C_0^\infty(\R^n)$ is chosen such that 
\begin{equation}\notag
\chi(\xi)=
\begin{cases}
1, \ \ \mbox { for } |\xi|\leq\frac{3}{4}\\
0, \ \ \mbox { for } |\xi|\geq 1.
\end{cases}
\end{equation}
Let 
\begin{equation}\notag
\varphi(\xi)=\chi(\xi/2)-\chi(\xi), \qquad \ \
\varphi_q(\xi)=
\begin{cases}
\varphi(\lambda_q^{-1}\xi)  \ \ \ \mbox { for } q\geq 0,\\
\chi(\xi) \ \ \ \mbox { for } q=-1.
\end{cases}
\end{equation}
For a tempered distribution vector field $v$ on the torus $\mathbb{T}^n$ we consider the Littlewood-Paley projections 
\bg\label{eq:LPP}
v_q(x)=\sum_{k\in\mathbb{Z}^n}\hat v(k)\varphi_q(k)e^{ik\cdot x}, \ \ \ \ q\geq -1.
\ed
The following Littlewood-Paley decomposition
\bg\notag
v=\sum_{q=-1}^\infty v_q
\ed
holds in the distribution sense. Essentially the sequence of the smooth functions $\varphi_q$ forms a dyadic partition of the unit. To simplify the notation, we denote
\bg\notag
v_{\leq q}=\sum_{j=-1}^qv_j, \ \ \ \ \tilde v_q=v_{q-1}+v_q+v_{q+1}.
\ed
By the definition of $\varphi_q$, it is noticed that $\supp (\varphi_p)\cap\supp (\varphi_{p'})=\emptyset$ if $|p-p'|\geq 2$.

By the Littlewood-Paley projection we define the Besov spaces $B_{r,l}^{s}$ on the torus $\mathbb T^3$ for $s\in\R$ and $1\leq l,r\leq\infty$. Denote the norm
\begin{align}\notag
&\|f\|_{B_{r,l}^{s}}=\left[\sum_{q\geq -1}(\lambda_q^s\|f_q\|_r)^l\right]^{1/l}.
\end{align}
Then
\begin{equation}\notag
B_{r,l}^s(\mathbb T^n)=\left\{f\in\mathit S': \|f\|_{\dot B_{r,l}^s}<\infty\right\},
\end{equation}
where $\mathit S'$ denotes the space of all tempered distributions.

We will often use the following inequality for the dyadic blocks of the Littlewood-Paley decomposition (see \cite{L}):
\begin{Lemma}\label{le:bern}(Bernstein's inequality)
For all $\alpha\in\mathbb {N}^n$, $q\in\mathbb{Z}$, $1\leq p\leq\infty$ and for all tempered distributions $f\in\mathcal S'$, we have
\bg\label{Bern}
\left\|\frac{\partial^\alpha}{\partial x^\alpha}f_q\right\|_{p}\sim\lambda_q^{|\alpha|}\|f_q\|_{p}.
\ed
\end{Lemma}

\subsection{The existence of Leray-Hopf type of weak solution to the incompressible MHD system}

We recall the result on the existence of weak solutions for the MHD system by Duvaut and Lions \cite{DL}. 
\begin{Theorem} \label{thm:ex}
For any  $(u_0, b_0)\in (L^2)^2$ there exists a weak solution $(u,b)$ to (\ref{MHD})- (\ref{init}) satisfying 
\[
u, b \in L^\infty(0,\infty; L^2) \cap L^2(0,\infty;H^1).
\]
\end{Theorem}
Moreover, there exists a weak solution satisfying the energy inequality
\bg\label{energy}
\|u(t)\|_2^2+\|b(t)\|_2^2+2\int_0^t\left(\mu\|\nabla u(s)\|_2^2+\nu\|\nabla b(s)\|_2^2\right) \,ds\leq\|u_0\|_2^2+\|b_0\|_2^2,
\ed
for all $t>0$.

\bigskip

\section{Construction of initial data}
\label{sec:data}

In this section we construct initial data $u_0\in B_{r,\infty}^{\frac{3}{r}+\theta-3}$ for the NSE and $(u_{0b},b_0)\in B_{r,\infty}^{\frac{3}{r}+\theta-3}\times B_{r,\infty}^{\frac{3}{r}+\gamma-3}$ for the MHD system with finite energy. The construction is similar to the one in \cite{CS}. 

Let  $\theta, \gamma>3/2$. We take any strictly decreasing sequence $\{q_j\}$ such that
\begin{equation}\label{qs}
\lambda_{q_i}^{4-\theta} \leq \lambda_{q_{i+1}}^{2\theta-3},
\ \ \qquad \lambda_{q_i}^{4-\theta}\leq \lambda_{q_{i+1}}^{2\gamma-3},
\ \ \qquad \lambda_{q_i}^{4-\gamma}\leq\lambda_{q_{i+1}}^{\theta+\gamma-3}.
\end{equation}
Given $c>0$, consider the following sets: 
\begin{align}\notag
&L_j=[\lambda_{q_j}, (1+c)\lambda_{q_j}]\times[-c\lambda_{q_j},c\lambda_{q_j}]^2\cap \mathbb{Z}^3\\
&M_j=[-c\lambda_{q_j-1},c\lambda_{q_j-1}]^2\times[\lambda_{q_j-1},(1+c)\lambda_{q_j-1}]\cap\mathbb{Z}^3\notag\\
&N_j=L_j+M_j\notag\\
&L_j^*=-L_j, \ \ M_j^*=-M_j, \ \ N_j^*=-N_j\notag.
\end{align}
Let $p(k)$ be the symbol of the Leray-Hopf projection
\bg\notag
p(k)=I-\frac{k\otimes k}{|k|^2}.
\ed
We denote
\begin{equation}\notag
\vec{\mathfrak e}_1(k)=p(k)\vec{e}_1, \ \ \  \vec{\mathfrak e}_2(k)=p(k)\vec{e}_2, \qquad k\in\mathbb{Z}^3\setminus \left\{0\right\},
\end{equation}
where $\vec e_j$ stands for a standard basis vector. Define
\begin{equation}\notag
\begin{split}
&\psi_1(k)=\vec{\mathfrak e}_2(k)\chi_{L_j\cup L_j^*}+i(\vec{\mathfrak e}_2(k)-\vec{\mathfrak e}_1(k))\chi_{N_j}-i(\vec{\mathfrak e}_2(k)-\vec{\mathfrak e}_1(k))\chi_{N_j^*},\\
&\psi_2(k)=\vec{\mathfrak e}_1(k)\chi_{M_j\cup M_j^*}. \\
\end{split}
\end{equation}

For the NSE, we choose initial data $u_0=U$, where 
\begin{equation}\label{u0}
\begin{split}
U=\sum_{j\geq 1}\lambda_{q_j}^{-\theta} \mathcal{F}\left(\psi_1(\xi)+ \psi_2(\xi) \right) 
\end{split}
\end{equation}
Due to the fact that $\phi(\xi)$ is flat around spheres $|\xi|=\lambda_q$, one can check
that for $c$ small enough we have
\[
\mathcal{F}(U_{q_j})(\xi)=\lambda_{q_j}^{-\theta}\psi_1(\xi), \qquad
\mathcal{F}(U_{q_j-1})(\xi)=\lambda_{q_j}^{-\theta}\psi_2(\xi), \qquad \mathcal{F}(U_{q_j+1})(\xi)=0.
\]
Hence $\tilde U_{q_j}=U_{q_j-1}+U_{q_j}$. It is also clear that
$\nabla\cdot u_0=0$.

For the MHD system, we choose initial data $u_{0b}$ and $b_0$ as
\begin{equation}\label{mhd-u0}
\begin{split}
u_{0b}=\sum_{j\geq 1}U_{q_j}, \qquad b_0=\sum_{j\geq 1}B_{q_j}, \qquad \mbox { with }\\
\mathcal{F}(B_{q_j})(\xi)=\lambda_{q_j}^{-\gamma}\psi_2(\xi).
\end{split}
\end{equation}
It also holds that
$$\nabla\cdot u_{0b}=\nabla\cdot b_0=0.$$


\begin{Lemma}\label{le:u0}
Let $\theta>3/2$.
For all $1< r\leq \infty$, we have $u_0, u_{0b}\in B_{r,\infty}^{\frac{3}{r}+\theta-3}$ and $b_0\in B_{r,\infty}^{\frac{3}{r}+\gamma-3}$. In particular, $u_0, u_{0b}\in H^{\theta-\frac32-s}$ and $b_0\in H^{\gamma-\frac32-s}$ for any $s>0$.
\end{Lemma}
\pf We only prove the conclusions for $u_0$. On the block $L_j$, for $1< r<\infty$, we have, by the boundedness of the Leray-Hopf projection and the $L^p$ estimate of Drichihlet kernel (see \cite{Gr08})
\begin{align}\notag
\|\lambda_{q_j}^{-\theta}\mathcal{F}^{-1}(\vec{e}_2(\xi)\chi_{L_j})\|_{r}&\lesssim \lambda_{q_j}^{-\theta}\|\mathcal{F}^{-1}(\chi_{L_j})\|_{r}\\
&\lesssim \lambda_{q_j}^{-\theta}\lambda_{q_j}^{3(1-\frac{1}{r})}.\notag
\end{align}
On all the other blocks, we have similar estimates. Hence,
\begin{equation}\notag
\lambda_{q_j}^{\frac{3}{r}+\theta-3}\|U_{q_j}\|_{r}\lesssim 1, \ \ \ \lambda_{q_j}^{\frac{3}{r}+\theta-3}\|U_{q_j-1}\|_{r}\lesssim 1.
\end{equation}
Therefore, $u_0\in \dot B_{r,\infty}^{\frac{3}{r}+\theta-3}$, for $1< r<\infty$. When $r=\infty$, 
\begin{align}\notag
\|U_{q_j}\|_{\infty}&\lesssim \|\mathcal {F}(U_{q_j})\|_{1}\\
&\lesssim \lambda_{q_j}^{-\theta}\left(\int_{\chi_{A_j\cup A_j^*}}1d\xi+\int_{\chi_{C_j}}1d\xi+\int_{\chi_{C_j^*}}1d\xi\right)\notag\\
&\lesssim \lambda_{q_j}^{3-\theta}\notag.
\end{align}
And similarly, we have 
\begin{equation}\notag
\|U_{q_j-1}\|_{\infty}\lesssim\lambda_{q_j}^{3-\theta}.
\end{equation} 
Therefore, $u_0\in \dot B_{\infty,\infty}^{\theta-3}$.
In particular, for $r=2$, the embedding $\dot B^{\theta-\frac{3}{2}}_{2,\infty}\subset H^{\theta-\frac32-s}$ holds for all $s>0$. Similar conclusion holds for $b_0$.
\cbdu

\begin{Remark}
Specifically the assumption $\theta,\gamma>\frac32$ implies that $u_0\in L^2$ and $(u_{0b},b_0)\in L^2\times L^2$ which indicates the initial data has finite energy.
\end{Remark}

As a consequence of Lemma \ref{le:u0}, one can see that, for $1<p\leq\infty$
\begin{equation}\label{eq:Ul2}
\|U_{q_j-1}\|_{p}+\|U_{q_j}\|_{p}\lesssim \lambda_{q_j}^{3-\theta-\frac{3}{p}},
\ \ \qquad \|B_{q_j}\|_{p}\lesssim \lambda_{q_j}^{3-\gamma-\frac{3}{p}}.
\end{equation}

The following estimates are essential to produce the discontinuity of the weak solutions.

\begin{Lemma}\label{le:tril}
Let $u_0$, $u_{0b}$, $b_0$ be defined as in (\ref{u0})-(\ref{mhd-u0}). Then the trilinear terms satisfy
\begin{equation}\notag
\begin{split}
\mathcal B(u_0, u_0, U_{q_j})\sim B(u_{0b}, u_{0b}, U_{q_j})\sim \lambda_{q_j}^{7-3\theta}, 
\ \ \qquad \mathcal B(b_0, b_0, U_{q_j})\sim \lambda_{q_j}^{7-2\gamma-\theta},\\
\mathcal B(u_{0b}, b_0, B_{q_j})\sim \lambda_{q_j}^{7-2\gamma-\theta}, 
\ \ \qquad \mathcal B(b_0, u_{0b}, B_{q_j})\sim \lambda_{q_j}^{7-2\gamma-\theta}.
\end{split}
\end{equation}
\end{Lemma}
\pf We only give a proof for the first one. The other estimates can be obtained in a similar way.
Note that $\supp \left\{\hat {U_{i}}\right\}\cap\supp \left\{\hat {U_{j}}\right\}=\varnothing$ for any $|i-j|\geq 2$. 
We decompose the term as
\begin{equation}\notag
\begin{split}
\mathcal B(u_0, u_0, U_{q_j})=&\sum_{k\geq j+1}\mathcal B(\tilde{U}_{q_k},\tilde{U}_{q_k}, {U}_{q_j})+\mathcal B(\tilde{U}_{q_j},\tilde{U}_{q_j}, {U}_{q_j})\\
&+\mathcal B(U_{\leq q_{j-1}},\tilde{U}_{q_j}, {U}_{q_j})+\mathcal B(\tilde{U}_{q_j}, U_{\leq q_{j-1}}, {U}_{q_j})\\
\sim&\sum_{k\geq j+1}\mathcal B(\tilde{U}_{q_k},\tilde{U}_{q_k}, {U}_{q_j})
+\mathcal B(U_{q_j-1}, U_{q_j}, {U}_{q_j})\\
&-\mathcal B(U_{q_j}, U_{q_j}, {U}_{\leq q_{j-1}})
\equiv I+II-III
\end{split}
\end{equation}
where we used integration by parts and the divergence free property of $U_{q_j}$.\\
Applying Bernstein's inequality (\ref{Bern}) and (\ref{eq:Ul2}) yields, for $\theta>3/2$
\begin{equation}\notag
|I|\lesssim\lambda_{q_j}\|U_{q_j}\|_{\infty}\sum_{k\geq j+1}\|\tilde{U}_{q_k}\|_{2}^2\lesssim\lambda_{q_j}^{4-\theta}\sum_{k\geq j+1}\lambda_{q_k}^{3-2\theta}\lesssim \frac{\lambda_{q_j}^{4-\theta}}{\lambda_{q_{j+1}}^{2\theta-3}}\leq \epsilon.
\end{equation}
Similarly,
\begin{equation}\notag
|III|\lesssim\|U_{q_j}\|_{2}^{3-2\theta}\sum_{k\leq j-1}\lambda_{q_k}\|{U}_{q_k}\|_{\infty}\lesssim \frac{\lambda_{q_{j-1}}^{4-\theta}}{\lambda_{q_{j}}^{2\theta-3}}\leq \epsilon.
\end{equation}
Using (\ref{Bern}) and (\ref{eq:Ul2}), the term $II$ is estimated as
\begin{equation}\notag
|II|\lesssim\lambda _{q_j}\|U_{q_j}\|_{2}^2\|U_{q_j}\|_{\infty}\sim \lambda _{q_j}^{7-3\theta}.
\end{equation}
The conclusion follows immediately.
\cbdu

\bigskip

\section{Discontinuous weak solutions to the NSE}
\label{sec:NSE}

\def \l {\lambda}
\def \lq {\l_{q_j}}
\def \lqm { \l_{q_j - 1} }
\def \lqk { \l_{q_k} }

In this section, we investigate the Navier-Stokes equation
\begin{equation}\label{NSE}
\begin{split}
&u_t- \mu\triangle u +u\cdot\nabla u+\nabla p=0,\\
&\nabla \cdot u=0, 
\end{split}
\end{equation}
with initial data $u_0$ given by (\ref{u0}). By an analogous analysis as in \cite{CS}, we show that the weak solutions of (\ref{NSE}) are discontinuous at initial time in a large class of Besov spaces, as stated in Theorem \ref{thm-nse}.

Denote $E(t) = \int_0^t \|\nabla u\|_2^2 ds$. Multiplying (\ref{NSE}) by $\tilde u_{q_j}$ and integrating over the space yields
\begin{equation}\label{crucial}
\begin{split}
\|\tilde{u}_{q_j}(t)\|_2^2 \geq &\|U_{q_j}\|_2^2 - \mu E(t) + c_1 \lq^{7-3\theta} t \\
&- c_2\int_0^t \left|  \mathcal B(u, u, u_{q_j}) - \mathcal B(U, U, U_{q_j})  \right| ds,
\end{split}
\end{equation}
for some positive constants $c_1$ and $c_2$. One contradiction argument will lead the conclusion of the theorem.

Suppose that for every $\delta >0$ there exists $t_0 = t_0(\delta) >0$ such that $\|u(t)- U\|_{B^{\frac 3r-3+\theta}_{r,\infty}} < \d$ for all $0<t\leq t_0$.
Denoting $w = u-U$, it follows
\begin{equation}\notag
\|w_p\|_r\leq\delta\lambda_p^{3-\theta-\frac 3r}, \ \ \mbox { for all } p\geq -1.
\end{equation}
After writing
\begin{equation}\notag
\begin{split}
&\mathcal B(u, u, u_{q_j}) - \mathcal B(U, U, U_{q_j})\\
= &\mathcal B(w, U, U_{q_j}) + \mathcal B(u, w, U_{q_j}) + \mathcal B(u, u, w_{q_j}) = A + B+ C,
\end{split}
\end{equation}
we estimate each term through the Bony's para-product (c.f. \cite{B}) decomposition as follows.
\begin{equation}\notag
\begin{split}
A = &\sum_{\substack{p',p'' \geq q_j \\ |p' - p''| \leq 2}} \mathcal B(w_{p'}, U_{p''}, U_{q_j}) + \mathcal B(w_{\leq q_j}, \tilde{U}_{q_{j}}, U_{q_j}) \\
&+ \mathcal B(\tilde{w}_{q_j}, U_{\leq q_j}, U_{q_j}) -r_A 
= A_1 + A_2 + A_3-r_A.
\end{split}
\end{equation}
with $r_A$ being the overlap of $A_2$ and $A_3$. Later $r_B, r_C, r_D, r_E, r_F$ have the same meaning.
Combining H\"{o}lder's inequality, and Bernstein's inequalities,  we obtain 
\begin{equation}\notag
\begin{split}
| A_1 | & \leq \|\nabla U_{q_j}\|_{\infty} \sum \|w_{p'}\|_{r} \|U_{p''}\|_{\frac r{r-1}} \lesssim \lq^{4-\theta} \sum \delta \l_{p''}^{3-2\theta}\lesssim  \d \lq^{7-3\theta}, 
\end{split}
\end{equation}
for $\theta>\frac32$;
\begin{equation}\notag
\begin{split}
|A_2| & = |\mathcal B(U_{q_j}, \tilde{U}_{q_{j}}, w_{\leq q_j} )| \leq \|U_{q_j}\|_{\infty}\|\tilde{U}_{q_j}\|_{\frac{r}{r-1}} \|\nabla w_{\leq q_j}\|_r \\
&\lesssim\lq^{3-2\theta+\frac3r} \sum_{p\leq q_j} \d\l_p^{4-\theta-\frac3r} 
\lesssim \d \lq^{7-3\theta}
\end{split}
\end{equation}
for $\theta+\frac3r<4$;
\begin{equation}\notag
\begin{split}
|A_3| & \leq \lq \|U_{q_j}\|_{\frac r{r-1}}\|U_{\leq q_j}\|_{\infty} \|\tilde{w}_{q_j}\|_{r}\\
&\lesssim \d \lq^{4-2\theta} \sum_{p\leq q_j}\l_p^{3-\theta} 
\lesssim \d \lq^{7-3\theta}
\end{split}
\end{equation}
for $\theta<3$.
We have shown
\begin{equation}\label{A}
|A| \lesssim \d \lq^{7-3\theta}, \ \quad \mbox { for } \frac32<\theta<3, \mbox { and } \theta+\frac3r<4.
\end{equation}
We decompose $B$ similarly,
\begin{equation}\notag
\begin{split}
B =& \sum_{\substack{p',p'' \geq q_j \\ |p' - p''| \leq 2}} \mathcal B(u_{p'}, w_{p''}, U_{q_j}) + \mathcal B(u_{\leq q_j}, \tilde{w}_{q_{j}}, U_{q_j}) \\
&+ \mathcal B(\tilde{u}_{q_j}, w_{\leq q_j},U_{q_j}) - r_B = B_1+B_2+B_3-r_B;
\end{split}
\end{equation}
\begin{equation}\notag
\begin{split}
\|B_1\| & \lesssim \lq \|U_{q_j}\|_{\frac{2r}{r-2}} \sum \|u_{p'}\|_2 \|w_{p''}\|_r
 \lesssim \d \lq^{\frac52-\theta+\frac3r} \sum_{p\geq\lq}\lambda_{p}^{2-\theta-\frac3r}\|\nabla u_{p}\|_2\\
& \lesssim \d \lq^{\frac92-2\theta} \sum_{p\geq\lq}\left(\frac{\lq}{\l_p}\right)^{\theta+\frac3r-2}\|\nabla u_{p}\|_2\lesssim \d \lq^{\frac92-2\theta}\|\nabla u_{p}\|_2
\end{split}
\end{equation}
for $\theta+\frac3r\geq 2$;
\begin{equation}\notag
|B_2|  =  \left| \mathcal B(U_{q_j}, \tilde{w}_{q_{j}}, u_{\leq q_j}) \right| \leq  \|U_{q_j}\|_{\frac{2r}{r-2}} \|\tilde{w}_{q_{j}}\|_r \|\nabla u_{\leq q_j}\|_2 
     \leq \d\lq^{\frac92-2\theta} \|\nabla u\|_2;
\end{equation}
\begin{equation}\notag
\begin{split}
|B_3| & \leq \|\tilde{u}_{q_j}\|_2 \|w_{\leq q_j}\|_r \|\nabla U_{q_j}\|_{\frac{2r}{r-2}} 
\lesssim \d\lq^{\frac32-\theta+\frac3r} \|\nabla\tilde{u}_{q_j}\|_2  \sum_{p \leq q_j}  \l_p^{3-\theta-\frac3r}  \\
&\lesssim 
\begin{cases}
\d\lq^{\frac92-2\theta} \|\nabla u\|_2,  \ \quad \mbox { if } \theta+\frac3r<3,\\
\d\lq^{\frac32-\theta+\frac3r} \|\nabla u\|_2, \ \quad \mbox { if } \theta+\frac3r>3,\\
\d\lq^{\frac32-\theta+\frac3r} q_j\|\nabla u\|_2, \ \quad \mbox { if } \theta+\frac3r=3.
\end{cases}
\end{split}
\end{equation}

We thus obtain
\begin{equation}\label{B}
|B| \lesssim 
\begin{cases}
\d\lq^{\frac92-2\theta} \|\nabla u\|_2,  \ \quad \mbox { if } 2\leq\theta+\frac3r<3,\\
\d\lq^{\frac32-\theta+\frac3r} \|\nabla u\|_2, \ \quad \mbox { if } \theta+\frac3r>3,\\
\d\lq^{\frac32-\theta+\frac3r} q_j\|\nabla u\|_2, \ \quad \mbox { if } \theta+\frac3r=3.
\end{cases}
\end{equation}
Similarly, 
\begin{equation}\notag
\begin{split}
C = &\sum_{\substack{p',p'' \geq q_j \\ |p' - p''| \leq 2}} \mathcal B(u_{p'}, u_{p''}, w_{q_j}) + \mathcal B(u_{\leq q_j}, \tilde{u}_{q_{j}}, w_{q_j}) \\
&+ \mathcal B(\tilde{u}_{q_j}, u_{\leq q_j}, w_{q_j}) -r_C = C_1+C_2+C_3-r_C;
\end{split}
\end{equation}
\begin{equation}\notag
\begin{split}
|C_1| & \leq \|\nabla w_{q_j}\|_r \sum_{p \geq q_j - 2} \|\tilde{u}_p\|_2\|\tilde{u}_p\|_{\frac{2r}{r-2}}
 \lesssim \d\lq^{4-\theta-\frac3r}  \sum_{p \geq q_j - 2}\l_p^{\frac3r-2} \|\nabla u_p\|_2^2\\
 & \lesssim \d \lq^{2-\theta} \|\nabla u\|_2^2, \ \quad \mbox { for } r\geq \frac32;\\
|C_2| & \leq \| \nabla u\|_{2} \|\tilde{u}_{q_{j}}\|_{\frac{2r}{r-2}} \|w_{q_j}\|_r \lesssim \d\lq^{2-\theta} \|\nabla u\|_2^2 ;\\
|C_3| &\lesssim \lq \|w_{q_j}\|_r \|\tilde{u}_{q_j}\|_{\frac{2r}{r-2}}  \|u\|_2\leq \d \lq^{3-\theta}  \|\nabla \tilde{u}_{q_j}\|_2.
\end{split}
\end{equation}
Thus,
\begin{equation}\label{C}
|C| \lesssim \d \lq^{2-\theta} \|\nabla u\|_2^2 + \d \lq^{3-\theta}  \|\nabla \tilde{u}_{q_j}\|_2.
\end{equation}

Combining \eqref{A}, \eqref{B}, \eqref{C} yields that, if $2\leq\theta+\frac3r<3$
\begin{equation} \notag
\begin{split}
&\int_0^{t_0} \left|  \mathcal B(u, u, u_{q_j}) - \mathcal B(U, U, U_{q_j})  \right| ds \\
\lesssim &\d \lq^{7-3\theta} t_0 + \d \lq^{9/2-2\theta} t^{1/2}_0 
+ \d\lq^{2-\theta} +  \d \lq^{3-\theta} \int_0^{t_0}\|\nabla \tilde{u}_{q_j}(s)\|_2\,ds\\
\lesssim &\d \lq^{7-3\theta} \left(t_0 + \lq^{\theta-5/2}t^{1/2}_0+\lq^{2\theta-5}+\lq^{2\theta-4}\int_0^{t_0}\|\nabla \tilde{u}_{q_j}(s)\|_2\,ds\right);
\end{split}
\end{equation} 
or if $3<\theta+\frac3r<4$, 
\begin{equation}\notag
\begin{split}
&\int_0^{t_0} \left|  \mathcal B(u, u, u_{q_j}) - \mathcal B(U, U, U_{q_j})  \right| ds \\
\lesssim &\d \lq^{7-3\theta} t_0 + \d \lq^{3/2-\theta+\frac3r} t^{1/2}_0 
+ \d\lq^{2-\theta} +  \d \lq^{3-\theta} \int_0^{t_0}\|\nabla \tilde{u}_{q_j}(s)\|_2\,ds\\
\lesssim &\d \lq^{7-3\theta} \left(t_0 + \lq^{2\theta+\frac3r-\frac{11}2}t^{1/2}_0+\lq^{2\theta-5}+\lq^{2\theta-4}\int_0^{t_0}\|\nabla \tilde{u}_{q_j}(s)\|_2\,ds\right);
\end{split}
\end{equation}
otherwise, if $\theta+\frac3r=3$,
\begin{equation}\notag
\begin{split}
&\int_0^{t_0} \left|  \mathcal B(u, u, u_{q_j}) - \mathcal B(U, U, U_{q_j})  \right| ds \\
\lesssim &\d \lq^{7-3\theta} \left(t_0 + \lq^{2\theta+\frac3r-\frac{11}2}q_jt^{1/2}_0+\lq^{2\theta-5}+\lq^{2\theta-4}\int_0^{t_0}\|\nabla \tilde{u}_{q_j}(s)\|_2\,ds\right).
\end{split}
\end{equation}
In the first case, we assume $\theta\leq 2$. Using the fact that
\[\int_0^{t_0} \|\nabla \tilde{u}_{q_j}(s)\|_2 ds \to 0, \quad \mbox { as } j \to \infty,\]
we can chose small enough $\delta$ and large enough $j_0$ such that for $j\geq j_0$
$$\int_0^{t_0} \left|  \mathcal B(u, u, u_{q_j}) - \mathcal B(U, U, U_{q_j})  \right| ds\leq\frac{c_1}{2c_2} \lq^{7-3\theta} t_0.$$
In the second and third cases, we assume 
\[\theta\leq 2, \quad 2\theta+\frac3r-\frac{11}2<0.\] Then for small enough $\delta$ and large enough $j_0$, the same estimate holds.
Going back to \eqref{crucial} it implies
$$
\|\tilde{u}_{q_j}(t_0)\|_2^2 \geq \|U_{q_j}\|_2^2 - \nu E(t_0) + c_1 \lq^{7-3\theta} t_0/2,
$$
for all $j >j_0$, which shows that $u(t_0)$ has infinite energy, a contradiction. In the end, we collect the conditions on the parameters and obtain that, either
\begin{equation}\notag
\begin{split}
&\frac32<\theta\leq 2, \ \qquad  2\leq\theta+\frac3r<3; \ \ \quad \mbox { or }\\
&\frac32<\theta\leq 2, \ \qquad r\geq\frac32,  \ \qquad  3\leq\theta+\frac3r<4,  \ \qquad  2\theta+\frac3r< \frac{11}2.
\end{split}
\end{equation}


\bigskip

\section{Discontinuous weak solutions to the MHD system}
\label{sec:MHD}

In the section, we show that, starting from data $(u_{0b},b_0)$ as defined in (\ref{mhd-u0}), a weak solution of the MHD system (\ref{MHD}) does not come back to $(u_0,b_0)$ in a large class of spaces. Namely we prove Theorem \ref{thm}.
 
Multiplying the second equation in (\ref{MHD}) by $b_{q_j}$ and integrating over the space and time interval $[0,t]$ yields
\begin{equation}\notag
\begin{split}
&\frac{1}{2}\|\tilde b_{q_j}(t)\|_{2}^2-\frac{1}{2}\|B_{q_j}\|_{2}^2\geq \int_{\mathbb{T}^3}\partial_t\tilde b_{q_j}\cdot b_{q_j}dx\\
&\geq-\nu\int_0^t\|\nabla b\|_{2}^2ds+\int_0^t\mathcal B(u, b, b_{q_j})-\mathcal B(b, u, b_{q_j})ds.
\end{split}
\end{equation}
Denote $E_b(t)=\int_0^t\nu\|\nabla b\|_{2}^2ds$. By Lemma \ref{le:tril} we have
\begin{equation}\label{ineq:energy}
\begin{split}
&\frac{1}{2}\|\tilde b_{q_j}(t)\|_{2}^2\geq\frac12\|B_{q_j}\|_{2}^2-E_b(t)+c_1\lambda_{q_j}^{7-2\gamma-\theta}t\\
&-c_2\int_0^t|\mathcal B(u, b, b_{q_j})-\mathcal B(u_{0b}, b_0, B_{q_j})|+|\mathcal B(b, u, b_{q_j})-\mathcal B(b_0, u_{0b}, B_{q_j})|ds\\
&\equiv \frac12\|B_{q_j}\|_{2}^2-E_b(t)+c_1\lambda_{q_j}^{7-2\gamma-\theta}t-c_2R(t)
\end{split}
\end{equation}
where $R$ represents the remainder term, and $c_1, c_2$ are positive constants. Again, we use a contradiction argument to show the conclusion of the theorem. 

Assume for every $\delta>0$ there exists $t_0=t_0(\delta)>0$ such that
\begin{equation}\label{ass}
\|u(t)-u_{0b}\|_{\dot B_{r,\infty}^{\theta-3+\frac 3r}}+\|b(t)-b_0\|_{\dot B_{r,\infty}^{\gamma-3+\frac 3r}}<\delta, \ \ \mbox { for all } 0<t\leq t_0.
\end{equation}
We claim that for a large enough $j_0$, the remainder term $R$ is bounded at $t_0$ as
\begin{equation}\label{claim}
R(t_0)<\frac{c_1}{2c_2}\lambda_{q_j}^{7-\theta-2\gamma}t_0, \ \ \mbox { for all } j\geq j_0.
\end{equation}

In the following we compute the second term in $R$ to obtain the desired estimate. The first term can be estimated similarly. 

Let $w=u-u_{0b}$ and $y=b-b_0$. Note that
by the assumption (\ref{ass}), we have
\begin{equation}\label{winfty}
\|w_p\|_r\leq\delta\lambda_p^{3-\theta-\frac 3r}, \ \ \ \|y_p\|_r\leq\delta\lambda_p^{3-\gamma-\frac 3r}, \ \ \mbox { for all } p\geq -1.
\end{equation}
The difference of the trilinear terms can be rewritten as
\begin{equation}\label{r4}
\begin{split}
&\mathcal B(b, u, b_{q_j})-\mathcal B(b_0, u_{0b}, B_{q_j})\\
=&\mathcal B(y, u_{0b}, B_{q_j})+\mathcal B(b, w, B_{q_j})+\mathcal B(b, u, y_{q_j})\notag\\
\equiv &D+E+F.\notag
\end{split}
\end{equation}
We decompose $D$ as
\begin{equation}\notag
\begin{split}
D=&\sum_{p',p''\geq {q_j},|p'-p''|\leq 2}\mathcal B(y_{p'}, (u_{0b})_{p''}, B_{q_j})+\mathcal B(y_{\leq {q_j}},  (\tilde u_{0b})_{q_j}, B_{q_j})\\
&+\mathcal B(\tilde y_{q_j}, (u_{0b})_{\leq {q_j}}, B_{q_j})-r_D\\
\equiv& D_1+D_2+D_3-r_D
\end{split}
\end{equation}
with $r_D$ being the overlap of $D_2$ and $D_3$. 

Applying the H\"older's inequality, (\ref{winfty}) and (\ref{eq:Ul2}) yields, 
\begin{equation}\notag
\begin{split}
|D_1|&\leq \|\nabla B_{q_j}\|_{\infty}\sum_{p',p''\geq {q_j},|p'-p''|\leq 2}\|y_{p'}\|_{r}\| (u_{0b})_{p''}\|_{\frac r{r-1}}\\
&\lesssim\lambda _{q_j}^{4-\gamma}\sum_{p',p''\geq {q_j},|p'-p''|\leq 2}\delta\lambda_{p'}^{3-\gamma-\frac3r}\lambda_{p''}^{\frac3r-\theta}\\
&\lesssim \delta\lambda_{q_j}^{7-2\gamma-\theta},
\end{split}
\end{equation}
for $\theta+\gamma>3$;
\begin{equation}\notag
\begin{split}
|D_2|&=|\mathcal B(B_{q_j}, (\tilde U)_{q_j}, y_{\leq{q_j}})|\leq\|B_{q_j}\|_{\infty}\|\tilde U_{q_j}\|_{\frac r{r-1}}\|\nabla y_{\leq{q_j}}\|_r\\
&\lesssim\lambda_{q_j}^{3-\theta-\gamma+\frac3r}\sum_{p\leq{q_j}}\delta\lambda_p^{4-\gamma-\frac3r}
\lesssim\delta\lambda_{q_j}^{7-2\gamma-\theta}
\end{split}
\end{equation}
for $\gamma+\frac3r<4$;
\begin{equation}\notag
\begin{split}
|D_3|&\leq \|\tilde y_{q_j}\|_{r}\|\nabla B_{q_j}\|_{\frac r{r-1}}\|U_{\leq{q_j}}\|_{\infty}
\lesssim\delta\lambda_{q_j}^{4-2\gamma}\sum_{p\leq {q_j}}\lambda_p^{3-\theta}
\lesssim\delta\lambda_{q_j}^{7-2\gamma-\theta}
\end{split}
\end{equation}
for $\theta<3$.
Hence,
\begin{equation}\label{est:d}
|D|\lesssim \delta \lambda_{q_j}^{7-2\gamma-\theta}.
\end{equation}

We decompose $E$ as,
\begin{align}\notag
E=&\sum_{p',p''\geq {q_j},|p'-p''|\leq 2}\mathcal B(b_{p'}, w_{p''}, B_{q_j})+\mathcal B(b_{\leq {q_j}}, \tilde w_{q_j}, B_{q_j})\\
&+\mathcal B(\tilde b_{q_j}, w_{\leq {q_j}}, B_{q_j})-r_E\notag\\
\equiv &E_1+E_2+E_3-r_E\notag
\end{align}
with $r_E$ being the overlap of $E_2$ and $E_3$. 

Similarly, we have, 
\begin{equation}\notag
\begin{split}
|E_1|&\leq \|\nabla B_{q_j}\|_{\frac{2r}{r-2}}\sum_{p',p''\geq {q_j},|p'-p''|\leq 2}\|b_{p'}\|_{2}\| w_{p''}\|_{r}\\
&\lesssim\lambda _{q_j}^{\frac52-\gamma+\frac 3r}\sum_{p'\geq {q_j}}\delta\lambda _{p'}^{2-\theta-\frac 3r}\|\nabla b_{p'}\|_{2}\\
&\lesssim\lambda _{q_j}^{\frac92-\theta-\gamma}\sum_{p'\geq {q_j}}\delta\left(\frac{\lambda _{q_j}}{\lambda _{p'}}\right)^{\theta+\frac 3r-2}\|\nabla b_{p'}\|_{2}\\
&\lesssim \delta\lambda_{q_j}^{\frac 92-\theta-\gamma}\|\nabla b\|_2
\end{split}
\end{equation}
for $\theta+\frac 3r\geq 2$;
\begin{equation}\notag
\begin{split}
|E_2|&=|\mathcal B(B_{q_j},\tilde w_{q_j}, b_{\leq{q_j}})|\leq\|B_{q_j}\|_{\frac{2r}{r-2}}\|\tilde w_{q_j}\|_r\|\nabla b_{\leq{q_j}}\|_2\lesssim \delta\lambda_{q_j}^{\frac{9}{2}-\theta-\gamma}\|\nabla b\|_2;
\end{split}
\end{equation}
\begin{equation}\notag
\begin{split}
|E_3|&\leq \|\tilde b_{q_j}\|_2\|w_{\leq{q_j}}\|_{r}\|\nabla B_{q_j}\|_{\frac{2r}{r-2}}\lesssim\|\nabla\tilde b_{q_j}\|_2\|B_{q_j}\|_{\frac{2r}{r-2}}\sum_{p\leq {q_j}}\delta\lambda_p^{3-\theta-\frac 3r}\\
&\lesssim \delta\lambda_{q_j}^{\frac{3}{2}-\gamma+\frac3r}\|\nabla b\|_2\sum_{p\leq {q_j}}\lambda_p^{3-\theta-\frac 3r}\\
&\lesssim 
\begin{cases}
\delta\lambda_{q_j}^{\frac{9}{2}-\theta-\gamma}\|\nabla b\|_2, \ \quad \mbox { if }  \theta+\frac 3r<3;\\
\delta\lambda_{q_j}^{\frac{3}{2}-\gamma+\frac3r}\|\nabla b\|_2, \ \quad \mbox { if }  \theta+\frac 3r> 3,\\
\delta\lambda_{q_j}^{\frac{3}{2}-\gamma+\frac3r}q_j\|\nabla b\|_2, \ \quad \mbox { if }  \theta+\frac 3r= 3.
\end{cases}
\end{split}
\end{equation}
Hence
\begin{equation}\label{est:e}
|E|\lesssim 
\begin{cases}
\delta\lambda_{q_j}^{\frac{9}{2}-\theta-\gamma}\|\nabla b\|_2, \ \quad \mbox { if } 2\leq\theta+\frac 3r<3;\\
\delta\lambda_{q_j}^{\frac{3}{2}-\gamma+\frac3r}\|\nabla b\|_2, \ \quad \mbox { if } \theta+\frac 3r> 3,\\
\delta\lambda_{q_j}^{\frac{3}{2}-\gamma+\frac3r}q_j\|\nabla b\|_2, \ \quad \mbox { if }  \theta+\frac 3r= 3.
\end{cases}
\end{equation}
An analogously decomposition for $F$ yields,
\begin{equation}\notag
\begin{split}
F=&\sum_{p',p''\geq {q_j},|p'-p''|\leq 2}\mathcal B(b_{p'}, u_{p''}, y_{q_j})+\mathcal B(b_{\leq {q_j}}, \tilde u_{q_j}, y_{q_j})\\
&+\mathcal B(\tilde b_{q_j}, u_{\leq {q_j}}, y_{q_j})-r_F\\
\equiv &F_1+F_2+F_3-r_F
\end{split}
\end{equation}
with $r_F$ being the overlap of $F_2$ and $F_3$. 

Again using the H\"older's inequality, Bernstein's inequalities (\ref{Bern}), (\ref{winfty}) and (\ref{eq:Ul2}) we infer that
\begin{equation}\notag
\begin{split}
|F_1|&\leq \|\nabla y_{q_j}\|_{r}\sum_{p',p''\geq {q_j},|p'-p''|\leq 2}\|b_{p'}\|_{2}\| u_{p''}\|_{\frac{2r}{r-2}}\\
&\lesssim\delta\lambda _{q_j}^{4-\gamma-\frac3r}\sum_{p',p''\geq {q_j},|p'-p''|\leq 2}\|\nabla b_{p'}\|_{2}\|\nabla u_{p''}\|_{2}\lambda_{p'}^{-1}\lambda_{p''}^{\frac3r-1}\\
&\lesssim\delta\lambda _{q_j}^{2-\gamma}\|\nabla b\|_{2}\|\nabla u\|_{2}\lesssim \delta
\lambda _{q_j}^{2-\gamma}(\|\nabla b\|_2^2+\|\nabla u\|_{2}^2)
\end{split}
\end{equation}
for $\frac3r\leq2$;
\begin{equation}\notag
\begin{split}
|F_2|&=|\mathcal B(y_{q_j}, \tilde u_{q_j}, b_{\leq{q_j}})|\leq\|y_{q_j}\|_r\|\tilde u_{q_j}\|_{\frac{2r}{r-2}}\|\nabla b_{\leq{q_j}}\|_2\\
&\lesssim\delta\lambda _{q_j}^{2-\gamma}\|\nabla\tilde u_{q_j}\|_2\|\nabla b_{\leq{q_j}}\|_2\lesssim \delta\lambda _{q_j}^{2-\gamma}(\|\nabla u\|_2^2+\|\nabla b\|_2^2),
\end{split}
\end{equation}
\begin{equation}\notag
\begin{split}
|F_3|&\leq \|\tilde b_{q_j}\|_{\frac{2r}{r-2}}\|u_{\leq{q_j}}\|_{2}\|\nabla y_{q_j}\|_{r}
\lesssim\lambda_{q_j}^{4-\gamma}\|\tilde b_{q_j}\|_2\|u\|_{2}\lesssim\delta\lambda_{q_j}^{3-\gamma}\|\nabla\tilde b_{q_j}\|_2.
\end{split}
\end{equation}
Thus, we have, for $r\geq\frac32$
\begin{equation}\label{est:f}
|F|\lesssim \delta\lambda _{q_j}^{2-\gamma}(\|\nabla u\|_2^2+\|\nabla b\|_2^2)+\delta \lambda_{q_j}^{3-\gamma}\|\nabla \tilde b_{q_j}\|_2.
\end{equation}

Combining (\ref{r4})--(\ref{est:f}) and the estimate $u,b\in L^2(H^1)$ gives that, if $2\leq\theta+\frac 3r<3$ 
\begin{equation}\notag
\begin{split}
&\int_0^{t_0}|\mathcal B(b, u, b_{q_j})-\mathcal B(b_0, u_{0b}, U_{q_j})|ds\\
&\lesssim\delta\left(\lambda_{q_j}^{7-\theta-2\gamma}t_0+\lambda_{q_j}^{\frac92-\theta-\gamma}t_0^{1/2}
+\lambda_{q_j}^{2-\gamma}+\lambda_{q_j}^{3-\gamma}\int_0^{t_0}\|\nabla \tilde b_{q_j}(s)\|_2\,ds\right)\\
&\lesssim\delta\lambda_{q_j}^{7-\theta-2\gamma}\left(t_0+t_0^{1/2}+\int_0^{t_0}\|\nabla \tilde b_{q_j}(s)\|_2\,ds\right)
\end{split}
\end{equation}
providing that $\gamma\leq\frac 52$ and $\theta+\gamma\leq 4$. Otherwise, if $\theta+\frac 3r\geq3$, 
\begin{equation}\notag
\int_0^{t_0}|\mathcal B(b, u, b_{q_j})-\mathcal B(b_0, u_{0b}, U_{q_j})|ds
\lesssim\delta\lambda_{q_j}^{7-\theta-2\gamma}\left(t_0+t_0^{1/2}+\int_0^{t_0}\|\nabla \tilde b_{q_j}(s)\|_2\,ds\right)
\end{equation}
for $\theta+\gamma+\frac3r<\frac {11}2$ and $\theta+\gamma\leq 4$.

Therefore we choose $j_0$ large enough and $\delta$ small enough such that, for all $j\geq j_0$ 
\begin{equation}\notag
\int_0^{t_0}|\mathcal B(b, u, b_{q_j})-\mathcal B(b_0, u_{0b}, U_{q_j})|ds\leq \frac{c_1}{4c_2}\lambda_{q_j}^{7-\theta-2\gamma}t_0.
\end{equation}
The first term in the integral $R$ can be estimated analogously and satisfies
\begin{equation}\notag
\int_0^{t_0}|\mathcal B(u, b, b_{q_j})-\mathcal B(u_{0b}, b_0, U_{q_j})|ds\leq \frac{c_1}{4c_2}\lambda_{q_j}^{7-\theta-2\gamma}t_0,
\end{equation}
for all $j\geq j_0$. Therefore, we have shown that the claim (\ref{claim}) holds under the assumption (\ref{ass}).
It follows from (\ref{ineq:energy}) and (\ref{claim}) that
\begin{equation}\notag
\frac{1}{2}\|\tilde b_{q_j}(t_0)\|_{2}^2\geq\frac12 \|B_{q_j}\|_{2}^2-E_b(t_0)+\frac{c_1}{2}\lambda_{q_j}^{7-\theta-2\gamma}t_0
\end{equation} 
 which implies $\|b(t_0)\|_2$ is infinity. It is a contradiction which is obtained under the conditions
 \[\frac32<\theta<3, \ \quad \frac32<\gamma\leq \frac52, \ \quad \theta+\gamma\leq 4, \ \quad 2\leq\theta+\frac3r<3, \ \ \quad \gamma+\frac3r<4;\]
 or
 \[r\geq\frac32, \, \frac32<\theta<3, \, \gamma>\frac32, \, \theta+\gamma\leq 4, \, \theta+\frac3r\geq3, \, \gamma+\frac3r<4, \, \theta+\gamma+\frac3r<\frac{11}2.\]
 
Anagolous analysis will give a contradiction that $\|u(t_0)\|_2$ is infinity at a certain time $t_0$ under an alternate assumption on the parameter triplet $(r,\theta,\gamma)$.

\bigskip



\end{document}